\documentclass[12pt]{amsart}

\newcommand{\C}{{\mathbb C}}
\newcommand{\Q}{{\mathbb Q}}
\newcommand{\Z}{{\mathbb Z}}
\newcommand{\R}{{\mathbb R}}
\newcommand{\A}{{\mathbb A}}
\newcommand{\D}{{\mathcal D}}
\newcommand{\Prin}{{\mathcal P}}

\newcommand{\tors}{_{\text{tors}}}
\newcommand{\GCD}{\operatorname{GCD}}
\newcommand{\Aut}{\operatorname{Aut}}
\newcommand{\End}{\operatorname{End}}
\newcommand{\ch}{\operatorname{char}}
\newcommand{\Gal}{\operatorname{Gal}}

\newcommand{\Lbar}{\overline{L}}

\newcommand{\Fq}{{{\mathbb F}_q}}
\newcommand{\Fqbar}{{\overline{\mathbb F}_q}}
\newcommand{\Ftwo}{{{\mathbb F}_2}}
\newcommand{\Fqm}{{\mathbb F}_{q^m}}
\newcommand{\Fqs}{{\mathbb F}_{q^s}}
\newcommand{\mm}{{\mathfrak m}}
\newcommand{\nn}{{\mathfrak n}}
\newcommand{\pp}{{\mathfrak p}}
\newcommand{\qq}{{\mathfrak q}}
\newcommand{\PP}{{\mathfrak P}}

\newcommand{\defequal}{\stackrel{\text{def}}{=}}
\newcommand{\nichts}{{\left.\right.}}
\newcommand{\isom}{\cong}

\newtheorem{lemma}{Lemma}
\newtheorem{cor}{Corollary}
\newtheorem{prop}{Proposition}
 
\newtheorem{theorem}{Theorem}

\theoremstyle{definition}
	
\newtheorem{conj}{Conjecture}

\theoremstyle{remark}
\newtheorem{rem}{Remark}

\begin{document}

\title[Torsion in Drinfeld Modules]{Torsion in Rank~1 Drinfeld Modules and the Uniform Boundedness Conjecture}
\author{Bjorn Poonen}
\thanks{This research was supported by an NSF Mathematical Sciences Postdoctoral Research Fellowship.  Research at MSRI is supported in part by NSF grant DMS-9022140.}
\address{Mathematical Sciences Research Institute \\ Berkeley, CA 94720-5070, USA}
\email{poonen@msri.org}
\date{June 19, 1995}

\begin{abstract}
It is conjectured that for fixed $A$, $r \ge 1$, and $d \ge 1$, there is a uniform bound on the size of the torsion submodule of a Drinfeld $A$-module of rank $r$ over a degree $d$ extension $L$ of the fraction field $K$ of $A$.
We verify the conjecture for $r=1$, and more generally for Drinfeld modules having potential good reduction at some prime above a specified prime of $K$.
Moreover, we show that within an $\Lbar$-isomorphism class, there are only finitely many Drinfeld modules up to isomorphism over $L$ which have nonzero torsion.
For the case $A=\Fq[T]$, $r=1$, and $L=\Fq(T)$, we give an explicit description of the possible torsion submodules.
We present three methods for proving these cases of the conjecture, and explain why they fail to prove the conjecture in general.
Finally, an application of the Mordell conjecture for characteristic $p$ function fields proves the uniform boundedness for the $\pp$-primary part of the torsion for rank~2 Drinfeld $\Fq[T]$-modules over a fixed function field.
\end{abstract}

\maketitle

\section{Conjectures and Theorems}
\label{theorems}

In a 1977 paper, Mazur~\cite{mazur} proved that if $E$ is an elliptic curve over $\Q$, its torsion subgroup is one of the following fifteen groups:
	$$\Z/N\Z, \qquad 1 \le N \le 10 \text{ or } N=12;$$
	$$\Z/2\Z \times \Z/2N\Z, \qquad 1 \le N \le 4.$$
In particular, there is a uniform bound for the size of the torsion subgroup of an elliptic curve over $\Q$.
(Although this in recent years has often been called Ogg's conjecture, it was essentially formulated by Levi~\cite{levi} in 1908, and was again formulated by Nagell~\cite{nagell} in 1949, before the appearance of Ogg's paper~\cite{ogg} in 1971.
See~\cite{sch} for more history.)

Recently, this has been generalized to arbitrary number fields.
Merel~\cite{merel} proved that for every $d \ge 4$, the order of the torsion subgroup of an elliptic curve $E$ over a number field $K$ of degree $d$ over $\Q$ is divisible only by primes less than $2^{d+1}(d!)^{5d/2}$.
Earlier results of Kamienny and Mazur~\cite{km} and Abramovich~\cite{abramovich} established bounds for $d \le 14$.
Also, Kamienny and Mazur~\cite{km} proved that the exponent of each such prime is bounded, with the bound depending only on $d$.
(Earlier work of Manin~\cite{manin} proved the boundedness of each exponent for fixed $K$.)
Together, these results imply that there is a bound $N(d)$ depending only on $d$ such that the order of the torsion subgroup of an elliptic curve over a number field of degree $d$ over $\Q$ is at most $N(d)$.
This paper considers the analogous questions for Drinfeld modules.

For an introduction to the theory of Drinfeld modules see~\cite{hayes} or Drinfeld's original paper~\cite{drinfeld}.
As usual in this theory, let $A$ denote the ring of regular functions on the affine curve obtained by deleting a closed point ``$\infty$'' from a nonsingular projective curve $X$ over $\Fq$.
Let $K$ be the fraction field of $A$.
In this paper, we will always equip a field extension $L$ of $K$ with the inclusion homomorphism $A \rightarrow L$.
In this context it makes sense to speak of rank $r$ Drinfeld $A$-modules over $L$, for each $r \ge 1$.
If $\phi$ is such a Drinfeld module, $^\phi L$ will denote the additive group of $L$ equipped with the $A$-module structure coming from $\phi$, and $(^\phi L)\tors$ will denote its torsion submodule; i.e.,
	$$(^\phi L)\tors = \{\, \alpha \in L \mid \phi_a(\alpha)=0 \text{ for some nonzero $a \in A$} \,\}.$$
Finally, $\#S$ will denote the cardinality of a set $S$.
It is known (see~\cite{denis2}, for example) that if $L$ is a finite extension of $K$, then $\# (^\phi L)\tors < \infty.$

The following two conjectures (a weak and a strong form) would give natural analogues of the recent results on uniform boundedness for elliptic curves.

\begin{conj}
\label{uniform}
For fixed $A$, $r \ge 1$, and finite extension $L$ of $K$, there is a uniform bound on $\#(^\phi L)\tors$ as $\phi$ ranges over rank $r$ Drinfeld $A$-modules over $L$.
\end{conj}

\begin{conj}
\label{uniform2}
For fixed $A$, $r \ge 1$, and $d \ge 1$, there is a uniform bound on $\#(^\phi L)\tors$ as $L$ ranges over extensions of $K$ of degree less than or equal to $d$, and $\phi$ ranges over rank $r$ Drinfeld $A$-modules over $L$.
\end{conj}

One might ask if there is a uniform bound over an even larger class of Drinfeld modules.
It is certainly not reasonable to ask for a uniform bound if $d$ is not fixed, because for any Drinfeld module $\phi$ over a function field $L$, the torsion submodule over $\Lbar$ is infinite.
If $L$ is fixed, but $r$ is allowed to vary, then again one cannot hope for a uniform bound on the torsion, because if $V$ is any finite sub-$\Fq$-vector space in $K=\Fq(T)$, then there is a polynomial of the form
	$$\phi_T(x) = Tx + a_1 x^q + a_2 x^{q^2} + \cdots a_r x^{q^r} \in K[x]$$
with zero set $V$, and this defines a Drinfeld $\Fq[T]$-module $\phi$ over $K$ whose torsion submodule is of size at least $\# V$.

The first main goal of this paper is to prove the following partial result.

\begin{theorem}
\label{rankone}
Conjecture~\ref{uniform2} holds when $r=1$ (for each $A$ and $d$).
\end{theorem}

In fact, we will show that for almost all rank~1 Drinfeld modules, there is no nontrivial torsion at all!

\begin{theorem}
\label{mostly}
For fixed $A$ and finite extension $L$ of $K$, there are at most finitely many rank~1 Drinfeld $A$-modules $\phi$ over $L$ (up to isomorphism over $L$) for which $(^\phi L)\tors \not= 0$.
\end{theorem}

Theorem~\ref{mostly} fails for rank $r \ge 2$.
For example, if $A=\Fq[T]$, $L=K=\Fq(T)$, then each $\alpha \in L^\ast$ gives rise to a rank~2 Drinfeld module
	$$\phi_T(x) = Tx - (\alpha+T) x^q + \alpha x^{q^2}$$
for which $1 \in \nichts^\phi L$ is torsion.
In fact, this gives infinitely many examples with good reduction at a specified place.

\medskip

For fixed $A$, there are only finitely many rank~1 Drinfeld modules over $L$ up to isomorphism over the algebraic closure $\Lbar$.
In fact, their number is the class number of $A$.
But each of these $\Lbar$-isomorphism classes may consist of infinitely many distinct classes of Drinfeld modules up to isomorphism {\em over $L$}, and these may have different torsion submodules, as in Theorem~\ref{specialcase} below, for example, so these theorems do not follow trivially even if $L$ is fixed.
On the other hand, the finiteness of the number of $\Lbar$-isomorphism classes implies that Theorems~\ref{rankone} and~\ref{mostly} are corollaries of the following two theorems, which apply to Drinfeld modules of arbitrary rank.
A {\em twist} of a Drinfeld module $\phi$ over $L'$ is a Drinfeld module over $L'$ which becomes isomorphic to $\phi$ over the algebraic closure.

\begin{theorem}
\label{general}
Let $\phi$ be a Drinfeld $A$-module over a finite extension $L$ of $K$.
Also fix $d \ge 1$.
There exists a uniform bound $C>0$ (depending on $\phi$, $A$, $L$, and $d$) such if $\psi$ is a twist of $\phi$ over a field extension $L'$ with $[L':L] \le d$, then $\# (^\psi L')\tors \le C$.
\end{theorem}

\begin{theorem}
\label{generalmostly}
If $\phi$ is a Drinfeld module over a finite extension $L$ of $K$, then at most finitely many of its twists $\psi$ over $L$ have $(^\psi L)\tors \not= 0$.
\end{theorem}

The elliptic curve analogue of Theorem~\ref{generalmostly}, namely that the torsion subgroups of almost all twists of a fixed elliptic curve over a number field are trivial, is not quite true, but becomes true if we ignore 2-torsion (and also 3-torsion in the case of zero $j$-invariant).

\medskip

In fact, we can generalize Theorem~\ref{general} further, to the following analogue of a known result for abelian varieties.
(To deduce Theorem~\ref{general} from the following, let $\pp$ be a prime of $K$ above which lies a prime of good reduction of $\phi$, and observe that every twist of $\phi$ then has potential good reduction at a place above $\pp$.)

\begin{theorem}
\label{potential}
Fix $r,d \ge 1$, fix $A$, and fix a prime $\pp$ of $A$.
There exists $C>0$ such that if $L$ is a finite extension of $K$ of degree at most $d$ and $\phi$ is a rank $r$ Drinfeld $A$-module over $L$ with potential good reduction at some prime above $\pp$, then $\#(^\phi L)\tors < C$.
\end{theorem}
But, as was remarked after the statement of Theorem~\ref{mostly}, we cannot generalize Theorem~\ref{generalmostly} in the same way.

\medskip

Sections~\ref{proofs},~\ref{proofs2} and~\ref{proofs3} give three different approaches to proving Theorem~\ref{rankone}.
The first two, which use local fields and adelic parallelotopes respectively, also prove Theorem~\ref{mostly} and the more general versions (Theorems~\ref{general} and~\ref{generalmostly}), and the first also proves Theorem~\ref{potential}.
The third, which uses the explicit description of the action of Galois on the abelian extensions obtained by adjoining torsion of sign-normalized rank~1 Drinfeld modules, proves only Theorem~\ref{rankone}.

Each method of proof has its own virtues.
(So there is reason for providing several proofs for the same theorems, beyond that of simply desiring to exhibit the various techniques that can be used to attack such problems.)
The proofs involving local fields are the quickest, and more importantly are the only ones which give Theorem~\ref{potential}.
The proofs involving adelic parallelotopes require only computations within the field over which the Drinfeld modules of interest are defined, and this makes them especially suitable for working out specific examples, as we will do in Section~\ref{examples}.
The proof involving class field theory gives the best bounds (just slightly worse than linear in the degree of the field, as opposed to exponentially large).
In fact, the bounds it gives will be proven to be best possible, up to a constant factor.

Section~\ref{higher} explains why these three methods fail to verify the conjectures for $r \ge 2$, and proves the following partial result for rank~2, which is analogous to Manin's local result for elliptic curves over number field mentioned above.

\begin{theorem}
\label{maninlike}
Let $A=\Fq[T]$ and $K=\Fq(T)$.
Fix a prime $\pp$ of $A$ and a finite extension $L$ of $K$.
Then there is a uniform bound on the size of the $\pp$-primary part of the torsion of rank~2 Drinfeld $A$-modules over $L$.
\end{theorem}

\section{Proofs using local fields}
\label{proofs}

As mentioned in the introduction, whenever we have a Drinfeld module over an extension of $K$, it is understood that the ring homomorphism $A \rightarrow L$ is the inclusion map.

\begin{prop}
\label{local}
Let $L$ be a finite extension of the completion $K_\pp$ of $K$ at a place $\pp$ other than $\infty$, and let $\phi$ be a Drinfeld module over $L$.
Then $(^\phi L)\tors$ is finite.
\end{prop}

\begin{proof}
If we fix a non-constant $T \in A$, then we may consider $\phi$ as a Drinfeld $\Fq[T]$-module over $L$, and the torsion submodule will be the same subset of $L$, since the torsion submodule can be characterized as the set of elements which have finite orbit under iteration of $\phi_T$.
Thus without loss of generality assume $A=\Fq[T]$.

Let $v$ be the discrete valuation on $L$.
If $v(x)$ is sufficiently negative, then $v(\phi_T(x))$ is more negative, and is determined by the leading term of $\phi_T$.
From this, if $g \in A=\Fq[T]$, $v(\phi_g(x))$ will be determined by the leading term in $g$.
In particular (still assuming $v(x)$ is sufficiently negative), $\phi_g(x)$ will not be zero.
Hence torsion elements are bounded below in valuation, and their closure is compact.
On the other hand, there are no nonzero torsion elements in a neighborhood of zero, because the action of $\phi_T$ near zero is locally analytically conjugate via a formal $A$-module isomorphism to multiplication-by-$T$, and $|T| \le 1$ with respect to the absolute value on $L$, so in a sufficiently small ball around 0, the entire Drinfeld $A$-action is analytically conjugate to the standard $A$-action.
Thus the group $(^\phi L)\tors$ is discrete as well as compact, so it is finite.
\end{proof}

\begin{rem}
The lemma could sometimes fail if $L$ were allowed to be $K_\infty$ or a finite extension thereof, as can be seen from the analytic uniformization theory.
(The same thing happens for elliptic curves: an elliptic curve over a $p$-adic field has finite torsion, whereas an elliptic curve over $\R$ or $\C$ has infinite torsion.)
\end{rem}

Theorem~\ref{potential} follows immediately from the following quantitative version of Proposition~\ref{local} for Drinfeld modules of potential good reduction.

\begin{theorem}
\label{localpotential}
Fix $r,d \ge 1$, fix $A$, and fix a prime $\pp$ of $A$.
There exists $C>0$ such that if $L$ is a finite extension of $K_\pp$ of degree at most $d$ and $\phi$ is a rank $r$ Drinfeld $A$-module over $L$ with potential good reduction, then $\#(^\phi L)\tors < C$.
\end{theorem}

\begin{proof}
Pick $T \in \pp$.
Any Drinfeld $A$-module can be considered as a Drinfeld $\Fq[T]$-module (with different but still bounded $r$ and $d$), so we may assume $A=\Fq[T]$, and $K_\pp=\Fq((T))$.
This time we equip any finite extension $L$ of $K_\pp$ with the valuation $v$ normalized by $v(T)=1$.
If $\phi$ is a rank $r$ Drinfeld $\Fq[T]$-module over $L$ with potential good reduction, then over the field obtained by adjoining the $(q^r-1)$-th root of the leading coefficient of $\phi_T$ to $L$, it is isomorphic to a Drinfeld module with good reduction.
Hence again replacing $d$ by something larger, but bounded, we may assume $\phi$ has good reduction over $L$.
In other words,
	$$\phi_T(x) = Tx + \text{(higher powers of $x$)} \in O[x],$$
where $O$ is the ring of integers of $L$.
All torsion of $\phi$ over $L$ is contained in $O$.
If $x$ is in the ideal $T^2O$, then $v(\phi_T(x))=v(x)+1$, and then by induction we see that $v(\phi_{T^k}(x))=v(x)+k$, which implies that the orbit of $x$ under $\phi_T$ is infinite, and hence $x$ cannot be torsion.
Thus $(^\phi L)\tors$ injects (as a group) into $O/T^2O$, which is of size $q^{2[L:K_\pp]} \le q^{2d}$.
\end{proof}

\begin{proof}[First proof of Theorem~\ref{mostly}]
Again we may assume $A=\Fq[T]$.
Write
	$$\phi_T(x) = Tx + a_1 x^q + a_2 x^{q^2} + \cdots a_r x^{q^r}.$$
Every twist $\psi$ of $\phi$ over $L$ is of the form
	$$\psi_T(x) = Tx + u^{q-1} a_1 x^q + u^{q^2-1} a_2 x^{q^2} + \cdots u^{q^r-1} a_r x^{q^r},$$
for some $u \in \Lbar^\ast$.
In fact, the requirement that $\psi_T$ have coefficients in $L$ forces $u^n \in L$, where $n$ is the $\GCD$ of $q^i-1$ over all $i$ for which $a_i \not=0$.
It is easy to see that $n$ also equals $\# \Aut(\phi)$.

Let $L_v$ be the completion of $L$ at some finite place, and let $F$ be the field extension obtained by adjoining to $L_v$ all $n$-th roots of all elements.
Since $n$ is prime to $\ch L_v$, Kummer theory implies that $F$ is a finite extension of $L_v$.
Every twist $\psi$ of $\phi$ is isomorphic to $\phi$ over $F$, so we obtain an isomorphism
	$$^\psi F \stackrel{u^{-1}}{\longrightarrow} \nichts^\phi F.$$
If $\psi$ and $\psi'$ are two such Drinfeld modules, then the images of $^\psi L$ and $^{\psi'} L$ intersect nontrivially only if $u \in L$, in which case $\psi$ and $\psi'$ are isomorphic over $L$.
Hence the union of $(^\psi L)\tors \setminus \{0\}$, with $\psi$ ranging over the distinct twists of $\phi$ over $L$, injects as a set into $(^\phi F)\tors$, which is finite by Proposition~\ref{local}.
Hence there can be only finitely many such twists $\psi$ for which $(^\psi L)\tors$ is nontrivial.
\end{proof}

\section{Proofs using adelic parallelotopes}
\label{proofs2}

The torsion submodule of a Drinfeld module over a global field can be characterized as the set of elements where the (global) canonical height is zero, or equivalently where each canonical local height is zero.
(See~\cite{denis2} for the global,~\cite{poonen} for the local.)
The latter condition can be satisfied at a given place $v$ only when the corresponding absolute value is small, and thus we can construct an adelic parallelotope in which all torsion lies.
Making these estimates uniformly over varying Drinfeld modules leads to proofs of our theorems.

For each place $\pp$ of a global function field $L$, let $\|\;\|_\pp$ denote the associated absolute value, normalized so that for each $x$ in the valuation ring $O_\pp$ of $\pp$, $\|x\| = \#(O_\pp/x)$.
The first lemma below bounds the torsion of a Drinfeld module within a certain adelic parallelotope, and the next bounds the number of field elements within such a parallelotope.

\begin{lemma}
\label{adele}
Let $\phi$ be a Drinfeld $A$-module over $L$.
Then for each place $\pp$ of $L$, there exists $c_\pp > 0$ such that for any extension $|\;|$ of $\|\;\|_\pp$ to a finite extension $E$ of $L$, every $\alpha \in (^\phi E)\tors$ satisfies $|\alpha| \le c_\pp$.
Moreover one can take $c_\pp=1$ for all but finitely many $\pp$.
\end{lemma}

\begin{proof}
The case $E=L$ is essentially a restatement of~(4) in Proposition~4 and~(3) in Proposition~1 of~\cite{poonen}, together with the observation (following from~(3) in Proposition~4 of~\cite{poonen}) that if some canonical local height function is nonzero at $\alpha$, then $\alpha$ is not in the torsion submodule.
It is clear from the proofs there that the constant $c_\pp$ depends only on the absolute values of the coefficients of the polynomials defining the Drinfeld module, so a single constant can be chosen to work for all $E$.
It is important here that $|\;|$ is literally an extension of $\|\;\|_\pp$ rather than the corresponding normalized absolute value, which is a power of $|\;|$.
\end{proof}

\begin{lemma}
\label{mkdivisor}
For each place $\pp$ of a global function field $L$, let $x_\pp$ be a positive real number, such that $x_\pp=1$ for all but finitely many $\pp$.
Fix one such place $\qq$, and let $k_\qq$ denote its residue field.
Let $M$ be the number of $\alpha \in L$ such that $\|\alpha\|_\pp \le x_\pp$ for all $\pp$.
Then either $M=1$ (i.e., the only such $\alpha$ is $0$), or
	$$M \le \# k_\qq \prod_\pp x_\pp.$$
\end{lemma}

\begin{proof}
This is a special case of Lemma~8 in~\cite{artinwhaples}.
The values of the constants in that lemma for this special case come out of its proof and case~3 of the proof of Lemma~7 there.
\end{proof}

\begin{proof}[Proof of Theorem~\ref{general}]
Just as in the proofs in Section~\ref{proofs}, we may assume that $A=\Fq[T]$, at the expense of increasing the rank.
Write
	$$\phi_T(x) = Tx + a_1 x^q + a_2 x^{q^2} + \cdots a_r x^{q^r}.$$
Every twist $\psi$ over $L'$ isomorphic to $\phi$ over $\Lbar$ is of the form
	$$\psi_T(x) = Tx + u^{q-1} a_1 x^q + u^{q^2-1} a_2 x^{q^2} + \cdots u^{q^r-1} a_r x^{q^r},$$
for some $u \in \Lbar^\ast$.
The torsion submodule of $\psi$ over $L'$ consists of those $\alpha \in L'$ for which $u\alpha$ is a torsion element of $\phi$.
Let $E=L'(u)$.

Suppose $u \alpha$ is a torsion element of $\phi$.
For each place $\pp$ of $L$, pick a constant $c_\pp$ as in Lemma~\ref{adele} for $\phi$.
Let $\PP$ be a place of $E$ above the places $\pp'$ of $L'$ and $\pp$ of $L$, and let $n(\PP/\pp)$ denote local degree.
Then $\|\;\|_\PP^{n(\PP/\pp)^{-1}}$ is an absolute value on $E$ extending $\|\;\|_\pp$.
So by Lemma~\ref{adele},
\begin{eqnarray}
\nonumber	\|u\alpha\|_\PP^{n(\PP/\pp)^{-1}}	& \le &	c_\pp,	\\
\label{bound1}	\|\alpha\|_\PP			& \le & c_\pp^{n(\PP/\pp')n(\pp'/\pp)} \|u\|_\PP^{-1},
\end{eqnarray}
Now fix $\pp'$.
If one takes the product of~(\ref{bound1}) over all places $\PP$ of $E$ above $\pp'$, one obtains
\begin{eqnarray}
\nonumber	\|\alpha\|_{\pp'}^{[E:L']}	& \le & c_\pp^{[E:L'] n(\pp'/\pp)} \prod_{\PP | \pp'} \|u\|_\PP^{-1},	\\
\label{bound2}	\|\alpha\|_{\pp'}		& \le & c_\pp^{n(\pp'/\pp)} \prod_{\PP | \pp'} \|u\|_\PP^{-1/[E:L']}.
\end{eqnarray}
Fix any prime $\qq$ of $L$, and let $\qq'$ be an extension to $L'$.
By Lemma~\ref{mkdivisor}, the number of $\alpha \in L'$ satisfying~(\ref{bound2}) for all places $\pp'$ of $L'$ is at most
\begin{eqnarray*}
	\lefteqn{ \# k_{\qq'} \prod_{\pp'} \left( c_\pp^{n(\pp'/\pp)} \prod_{\PP | \pp'} \|u\|_\PP^{-1/[E:L']} \right) }	\\
		& \le &	\left(\# k_\qq \right)^{[L':L]} \left( \prod_\pp c_\pp^{[L':L]} \right) \left( \prod_\PP \|u\|_\PP^{-1/[E:L']} \right)	\\
		& = & \left(\# k_\qq \prod_\pp c_\pp \right)^{[L':L]},
\end{eqnarray*}
by the product formula in $E$ applied to $u$.
Thus the size of the torsion submodule of $\psi$ over $L'$ is at most $C^{[L':L]} \le C^d$, where
	$$C=\# k_\qq \prod_\pp c_\pp.$$
\end{proof}

\begin{proof}[Second proof of Theorem~\ref{mostly}]
Again we may assume $A=\Fq[T]$.
Write
	$$\phi_T(x) = Tx + a_1 x^q + a_2 x^{q^2} + \cdots a_r x^{q^r}.$$
Every twist $\psi$ of $\phi$ over $L$ is of the form
\begin{equation}
\label{form}
	\psi_T(x) = Tx + f^{(q-1)/n} a_1 x^q + f^{(q^2-1)/n} a_2 x^{q^2} + \cdots + f^{(q^r-1)/n} a_r x^{q^r},
\end{equation}
where $n$ is as in the first proof of Theorem~\ref{general} and $f \in L^\ast$, and the isomorphism type of $\psi$ depends only on the image of $f$ in $L^\ast/L^{\ast n}$.

We claim that there exists a finite set of places $S$ of $L$ and integers $m_v$ for each $v \in S$ such that any $g \in L^\ast$ can be multiplied by an $n$-power to obtain $f$ satisfying
\begin{enumerate}
	\item $v(f) \le m_v$ if $v \in S$.
	\item $0 \le v(f) \le n-1$ if $v \not\in S$.
\end{enumerate}

Let $\D$ denote the group of degree zero divisors on the curve of which $L$ is the function field, and let $\Prin$ denote the subgroup of principal divisors.
Then the class group $\D/\Prin$ of $L$ is finite, and multiplication by $n$ maps it surjectively onto $n\D/n\Prin$, so $n\D/n\Prin$ is finite as well.
Let $D_1,\ldots,D_s$ be representatives for the cosets of $n\Prin$ in $n\D$.
Let $S$ be the set of places of $L$ which occur in $D_1,\ldots,D_n$, and let $v_0$ be one of these places.
(If $S$ is empty, enlarge it so it contains at least one place.)
Given $g$, with principal divisor $(g)=\sum \alpha_v v$, let $D$ be the degree zero divisor $D = \sum \alpha_v' v$ where $\alpha_v' = n \lfloor \alpha_v/n \rfloor$ for $v' \not= v_0$, and $\alpha_{v_0}'$ is set to make the total degree 0.
Then by definition of the $D_i$, we have $-D=D_i+n(u)$ for some $i$ and some $u \in L^\ast$.
So $f \defequal gu^n$ satisfies
	$$(f)=(g)+n(u)=(g)-D'-D_i.$$
The coefficient $\beta_v$ of $v$ in $(f)$ is $\alpha_v \bmod n$ if $v \not\in S$.
If $v \in S$ but $v \not= v_0$, then $\beta_v$ is bounded above and below, where the bound depends on the $D_i$ but not on $g$.
Since $\sum_v \beta_v =0$, and we have lower bounds on each $\beta_v$ for $v \not=v_0$ which are zero for all but finitely many $v$, we get also an upper bound on $\beta_{v_0}$, as desired.

From now on, we only consider twists $\psi$ of the form~(\ref{form}) with $f$ satisfying the two conditions of the previous paragraph.
We claim there is an adelic parallelotope of $L$ (of finite volume) which contains all torsion elements of all such $\psi$.
Let $T$ be the (finite) set of places of $L$ which occur in the factorizations of the nonzero $a_i$.

First suppose $v \not\in S \cup T$.
Then $v(f) \le n-1$ by condition~1 above.
If $v(x) \le -1$, then
\begin{align*}
	v(f^{(q^r-1)/n} a_r x^{q^r}) - v(f^{(q^i-1)/n} a_i x^{q^i})
	&= (q^r-q^i)/n \cdot v(f) + (q^r-q^i) v(x)	\\
	&\le (q^r-q^i)(n-1)/n - (q^r-q^i)	\\
	&< 0,
\end{align*}
so the term $f^{(q^r-1)/n} a_r x^{q^r}$ is larger $v$-adically than all the other terms in~(\ref{form}).
In particular it is larger than the term $Tx$, so $v(\psi_T(x))<v(x)$.
Hence $x$ tends to infinity $v$-adically under iteration of $\psi_T$, and so cannot be a torsion element.
In other words, the torsion of $\psi$ over $L$ lies within the set where $v(x) \ge 0$.

Now suppose $v \in S \cup T$.
Since $v(f) \le m_v$, we still have
	$$v(f^{(q^r-1)/n} a_r x^{q^r}) - v(f^{(q^i-1)/n} a_i x^{q^i}) < 0$$
provided that $v(x)<c_v$ where $c_v$ is sufficiently negative.
Then the same argument in the previous paragraph shows that the torsion of $\psi$ lies within the set where $v(x) \ge c_v$.

Thus we have shown that all torsion of $\psi$ over $L$ lies within an adelic parallelotope whose bounds are independent of $f$, hence independent of $\psi$.
By Lemma~\ref{mkdivisor}, there are only finitely many elements $\alpha_1,\ldots,\alpha_M$ within this adelic parallelotope.
If $\alpha_i$ is a nonzero torsion element for $\psi$ for a particular $f$, then $\psi_T(\alpha_i)$ is also torsion, so $\psi_T(\alpha_i)=\alpha_j$ for some $j$.
Rewriting this using~(\ref{form}) shows that $f$ satisfies a algebraic equation, so there are finitely many possibilities for $f$, for each $\alpha_i$ and $\alpha_j$.
The (finite) union of all these over $i$ and $j$ lists every $f$ for which the twist $\psi$ might have a nonzero torsion element.
\end{proof}

\section{Proofs using explicit class field theory}
\label{proofs3}

The bounds on the size of the torsion submodule of rank~1 Drinfeld modules given by the proofs of the last two sections grow exponentially in the degree $d$.
In this section we will give a third proof of Theorem~\ref{rankone}, and this one will give a much better bound (just slightly worse than linear in $d$).
In fact, the bound we obtain here will be shown to be best possible, up to a constant factor.

The basic tool we use in this section is the explicit description of the action of Galois on the torsion of a sign-normalized Drinfeld module.
This lets us determine which torsion elements belong to a particular global field.
The main drawback of these methods is that they apply only to rank~1 Drinfeld modules, and hence do not also give a proof of Theorem~\ref{general} on twists of Drinfeld modules of higher rank.
The quantitative version of Theorem~\ref{rankone} we prove is the following:

\begin{theorem}
\label{strengthening}
For fixed $A$, there is a constant $C>0$ such that if $\phi$ is a rank~1 Drinfeld $A$-module over a global field $L$ with $[L:K]=d$, then
	$$\#(^\phi L)\tors \le C \cdot d \log\log d.$$
For each $A$, this bound is best possible up to a constant factor.
\end{theorem}

In the proof, we will need bounds on a characteristic $p$ analogue of the Euler phi-function.
If $\mm$ is an ideal of $A$, we let $|\mm| = \# A/\mm$, $\deg \mm = \dim_{\Fq} A/\mm = \log_q |\mm|$, and $\varphi_A(\mm) = \# (A/\mm)^\ast$.
Just as in the classical case, we have
\begin{equation}
\label{classical}
	\varphi_A(\mm) = |\mm| \cdot \prod_{\pp | \mm} (1-1/|\pp|).
\end{equation}

\begin{lemma}
\label{liminf}
For each $A$,
   $$\liminf_{|\mm| \rightarrow \infty} \frac{\varphi_A(\mm) \log \log |\mm|} {|\mm|}$$
is a (finite) positive real number.
\end{lemma}

\begin{proof}
First we estimate the number of primes of $A$ of given degree.
Recall that $A$ is the ring of regular functions on a nonsingular projective curve $X$ minus a closed point $\infty$.
Degree $m$ primes of $A$ correspond to $\Gal(\Fqbar/\Fq)$-conjugacy classes of points of $X$ (other than $\infty$) defined over $\Fqm$ but not over any smaller field.
By the Riemann hypothesis for curves over finite fields, $\# X(\Fqm) = q^m + O(q^{m/2})$ as $m$ goes to infinity.
(The constant implied by the big-$O$ depends on $X$.)
Of the points defined over $\Fqm$, the number defined over smaller fields is at most
	$$\sum_{s|m,s \not=m} \# X(\Fqs) \le \sum_{s=1}^{m/2} \left[ q^s + O(q^{s/2}) \right] = O(q^{m/2})$$
so there remain $q^m + O(q^{m/2})$ whose field of definition is $\Fqm$.
Hence the number of primes of $A$ of degree $m$ is $q^m/m + O(q^{m/2}/m)$.
(It does not matter that we ignored the single prime $\infty$.)

Consider the product $\mm$ of the primes of degree up to $n$.
As is standard, $\Theta(f(n))$ denotes a function which is bounded above and below by positive multiples of $f(n)$.
We have
\begin{gather}
\nonumber	\deg \mm	= \sum_{m=1}^n (q^m/m + O(q^{m/2}/m)) m = \Theta(q^n)	\\
\label{loglog}	\log \log |\mm|	= \log \log q^{\Theta(q^n)} = n \log q + O(1)
\end{gather}
and
\begin{align}
\nonumber	\sum_{\deg \pp \le n} \log (1-1/|\pp|)
				&= \sum_{m=1}^n (q^m/m + O(q^{m/2}/m)) \log (1-q^{-m}) \\
\nonumber			&= \sum_{m=1}^n (q^m/m + O(q^{m/2}/m)) (-q^{-m}+O(q^{-2m})) \\
\nonumber			&= \sum_{m=1}^n (-1/m + O(q^{-m/2})) \\
\nonumber			&= - \log n + O(1)	\\
\label{phiprod}	\prod_{\deg \pp \le n} (1-1/|\pp|)
			&= \Theta(1/n).
\end{align}
(The estimate~(\ref{phiprod}) is an analogue of Mertens' Theorem which was proved independently by Rosen~\cite{rosen}, who also computes the constant implied by the $\Theta$ notation.)
By~(\ref{classical}),~(\ref{loglog}) and~(\ref{phiprod}),
	$$\frac{\varphi_A(\mm) \log \log |\mm|} {|\mm|} = (n \log q + O(1)) \Theta(1/n) = \Theta(1),$$
which proves that the $\liminf$ is finite.

On the other hand, for any ideal $\mm$ of degree $n$, the number $k$ of primes $\pp_1,\ldots,\pp_k$ dividing $\mm$ having degree at least $\log_q n$ is at most $n/\log_q n$, so by~(\ref{classical}),
\begin{align*}
	\varphi_A(\mm) / |\mm|	&= \prod_{i=1}^k (1-1/|\pp_i|) \prod_{\pp | \mm, \deg \pp <\log_q n} (1-1/|\pp|)	\\
				&\ge (1-1/q^{\log_q n})^{n/\log_q n} \prod_{\deg \pp < \log_q n} (1-1/|\pp|)	\\
				&\ge (1-1/\log_q n) \Theta(1/\lfloor \log_q n \rfloor) \quad\text{(by~(\ref{phiprod}))}	\\
				&= \Theta(1/\log n) \\
				&= \Theta(1/\log \log_q |\mm|) \\
				&= \Theta(1/\log \log |\mm|),
\end{align*}
which proves that the $\liminf$ is positive.
\end{proof}

\begin{cor}
\label{phigrowth}
For each $A$, there is a constant $c>0$ such that
   $$|\mm| \le c \,\varphi_A(\mm) \log \log \varphi_A(\mm)$$
for every ideal $\mm$ of $A$.
\end{cor}

\begin{proof}[Proof of Theorem~\ref{strengthening}]
Fix a sign-function, and let $H_A^+$ be the corresponding normalizing field for rank~1 Drinfeld $A$-modules, which is the minimal field of definition for each sign-normalized Drinfeld $A$-module.
(See~\cite{hayes} for definitions.)
Suppose $\phi$ is a rank~1 Drinfeld $A$-module over $L$, and $[L:K]=d$.
By Theorem~12.3 in~\cite{hayes}, $\phi$ is isomorphic to some sign-normalized Drinfeld module $\psi$, which is defined over $H_A^+$.
Let $E$ be a compositum of $L$ and $H_A^+$ in some algebraic closure of $K$.
Then $\phi$ and $\psi$ are twists over $E$, and they become isomorphic over a field extension $F$ obtained by adjoining a single $n$-th root, where $n$ (as in the proofs of Section~\ref{proofs}) equals $\# \Aut(\phi) = q-1$, by Corollary~5.14 in~\cite{hayes}.
In particular,
	$$(^\phi L)\tors \subset (^\phi F)\tors \isom (^\psi F)\tors.$$

Since $\psi$ is of rank 1, $(^\psi F)\tors$ is isomorphic as an $A$-module to $A/\mm$ for some ideal $\mm$ of $A$.
Let $K_\mm$ be the narrow ray class field, which is the field generated over $H_A^+$ by the $\mm$-torsion of $\psi$.
Section~16 of~\cite{hayes} gives us an isomorphism from $\Gal(K_\mm/H_A^+)$ to $(A/\mm)^\ast$, so
\begin{align*}
	\varphi_A(\mm)	&= \# \Gal(K_\mm/H_A^+)	\\
			&= [K_\mm:H_A^+]	\\
			&\le [F:H_A^+]		\\
			&\le n [E:H_A^+] \quad\text{(by definition of $F$)}\\
			&= n [LH_A^+:KH_A^+] \\
			&\le nd.
\end{align*}
Finally
\begin{align*}
	\# (^\phi L)\tors	&\le \# (^\psi F)\tors	\\
				&= |\mm|	\\
				&\le c \varphi_A(\mm) \log\log \varphi_A(\mm) \quad\text{(by Corollary~\ref{phigrowth})}	\\
				&\le c (nd) \log\log(nd)	\\
				&\le C \cdot d \log\log d,
\end{align*}
for some constant $C$ depending only on $A$.

This is best possible, because for the sequence of $\mm$'s which make the $\liminf$ in Lemma~\ref{liminf} small, if $\phi$ is a sign-normalized Drinfeld $A$-module over the narrow ray class field $K_\mm$, then
	$$d=[K_\mm:K]=[K_\mm:H_A^+][H_A^+:K]=O(\varphi_A(\mm))=O(|\mm|/\log\log|\mm|)$$
while
	$$\# (^\phi L)\tors \ge |\mm|.$$
\end{proof}

\section{Examples}
\label{examples}

Theorem~\ref{mostly} tells us in particular that over $K=\Fq(T)$, there are only finitely many rank~1 Drinfeld $\Fq[T]$-modules that have nonzero torsion.
Here we list them explicitly.

\begin{theorem}
\label{specialcase}
Let $A=\Fq[T]$ and $K=\Fq(T)$.
Let $\phi$ be a rank~1 Drinfeld $A$-module over $K$.

\medskip
\noindent{\em Case 1: $q=2$}

In this case $\phi$ is isomorphic to the Carlitz module $\phi_T(x)=Tx+x^2$ over $K$.
For the Carlitz module, $(^\phi K)\tors = \Ftwo + \Ftwo \cdot T$, and the torsion is isomorphic to $\Ftwo[T]/(T^2+T)$ as an $\Ftwo[T]$-module.

\medskip
\noindent{\em Case 2: $q>2$}

If $\phi_T(x)=Tx-(T+c)x^q$ for some $c \in \Fq$, then the torsion submodule $(^\phi K)\tors$ equals $\Fq \subset \nichts^\phi K$ and is isomorphic to $\Fq[T]/(T+c)$ as an $\Fq[T]$-module.  If $\phi$ is not isomorphic over $K$ to such a Drinfeld module, then $(^\phi K)\tors = 0$.
\end{theorem}

We will specialize the second proof of Theorem~\ref{generalmostly}, since it is the only one that does not require that we work within field extensions larger than $K$.
The fact that $K$ has class number~1 will simplify matters greatly.

\begin{proof}
Let $v_g$ denote the $\Z$-valued discrete valuation on $K$ associated with an irreducible element $g$ of $\Fq[T]$.
The only nontrivial valuation not of this form is the one at $\infty$, the degree function $\deg$.

We have
	$$\phi_T(x) = T x - f x^q$$
for some $f \in \Fq(T)$.
Since $\phi$ is isomorphic over $K$ to any Drinfeld module $\psi$ with
	$$\psi_T(x) = T x - u^{q-1} f x^q$$
for $u \in K^\ast$, we may assume without loss of generality that $f \in \Fq[T]$.
Also, since $\Fq[T]$ is a UFD, we may assume that $v_g(f) < q-1$ for each irreducible $f \in \Fq[T]$.
(When $q=2$, this means that $f=1$ and $\phi$ is the Carlitz module.)
The set $S \cup T$ in the second proof of Theorem~\ref{mostly} for us is just $\{\infty\}$.
Then that proof gives $v_g(\alpha) \ge 0$ for all irreducible $g$ in $\Fq[T]$.
In other words, $(^\phi K)\tors \subset \Fq[T]$.

Suppose $\alpha \subset \Fq[T]$ is of degree at least~1 if $q \not=2$, and at least~2 if $q=2$.
Then $\deg \alpha^{q-1} \ge 2$, so
\begin{eqnarray*}
	\deg (f \alpha^q)	& > &	\deg(T \alpha)	\\
	\deg(\phi_T(\alpha))	& = &	\deg (f \alpha^q) > \deg \alpha	\\
	\deg(\phi_{T^n}(\alpha))& \rightarrow & \infty
\end{eqnarray*}
as $n$ goes to $\infty$.
This implies that $\alpha$ is not a torsion element, because a torsion element has finite orbit under iteration $\phi_T$.
Thus $(^\phi K)\tors \subset \Fq$ if $q \not=2$, and $(^\phi K)\tors \subset \Ftwo + \Ftwo \cdot T$ if $q=2$.

For $q=2$ (and $f=1$),
\begin{eqnarray*}
	\phi_{T+1}(1)	& = & T	\\
	\phi_T(T)	& = & 0
\end{eqnarray*}
so $(^\phi K)\tors = \Ftwo + \Ftwo \cdot T$ and it is a cyclic module generated by~1 and isomorphic to $\Ftwo[T]/(T^2+T)$.

For $q>2$, $(^\phi K)\tors \subset \Fq$.
If $f=T+c$ with $c \in \Fq$, then $\phi_{T+c}(1)=0$, so it is clear that $(^\phi K)\tors = \Fq$, and that it is isomorphic to $\Fq[T]/(T+c)$ as an $\Fq[T]$-module.
Otherwise $\phi_T(1) = T-f \not\in \Fq$, so $1 \not\in (^\phi K)\tors$, and $(^\phi K)\tors=0$.
\end{proof}

\section{Uniform bounds for Drinfeld modules of higher rank}
\label{higher}

In this section we will first explain why the methods of the previous three section fail to prove Conjectures~\ref{uniform} or~\ref{uniform2} in general for Drinfeld modules of rank $r \ge 2$, and then prove Theorem~\ref{maninlike}, which is the only result beyond Theorem~\ref{potential} that we have in this situation.
If we work over a local field (the completion of the global function field at a prime not above $\infty$), we still have finiteness of the torsion submodule for each individual Drinfeld module by Proposition~\ref{local}, but on the other hand we will be able to exploit the Tate uniformization Drinfeld modules in the next proposition to construct examples where the torsion submodule over the local field can be arbitrarily large.
Let $H_A$ denote the Hilbert class field of $A$, which is the maximal unramified abelian extension of $K$ in which $\infty$ splits completely.  (See Section~15 of~\cite{hayes}.)

\begin{prop}
Let $L$ be a finite extension of the completion of $H_A$ at a place not above $\infty$.
Then for each $r \ge 2$ and nonzero $a \in A$, there exists a rank~$r$ Drinfeld module $\phi$ over $L$ with a torsion element of order $a$.
\end{prop}

\begin{proof}
There exists a rank~1 Drinfeld $A$-module $\psi$ over $L$, since the minimal field of definition for any rank~1 Drinfeld $A$-module is $H_A$.  (See~\cite{hayes}.)
By replacing $\psi$ with an isomorphic Drinfeld module, we may assume the coefficients of $\psi_a$ for all $a$ lie in the valuation ring $O$ of $L$.
Then $O$ and $L$ both inherit $A$-module structures from $\phi$, and the quotient $^\phi L/\nichts^\phi O$ as an $A$-module is isomorphic to the direct sum of a free module of rank $\aleph_0$ with a finite torsion module, by Theorem~2 in~\cite{poonen}.
Let $\lambda_1,\lambda_2,\ldots,\lambda_{r-1}$ be elements of $^\phi L$ whose images in $^\phi L/\nichts^\phi O$ are $A$-independent.
The $A$-module $\Lambda$ generated by the $\lambda_i$ is discrete since its intersection with $^\phi O$ is trivial.
Let $\phi_a(\Lambda)$ denote the submodule generated by $\phi_a(\lambda_1), \phi_a(\lambda_2), \ldots, \phi_a(\lambda_{r-1})$.
Then the quotient of $\psi$ by $\phi_a(\Lambda)$ can be given the structure of a rank~$r$ Drinfeld module over $L$, as in Section~7 of~\cite{drinfeld} on Tate uniformization, and the image of $\lambda_1$ in this quotient has order $a$.
\end{proof}

The adelic parallelotope method fails when $r \ge 2$, simply because the volume of the parallelotope needed to contain the torsion can be arbitrarily large.
For example, suppose $A=\Fq[T]$, $L=K=\Fq(T)$, $V$ is the sub-$\Fq$-vector space of $K$ with basis $\lambda_1,\ldots,\lambda_r$, and $\phi$ is the rank~$r$ Drinfeld $A$-module with $\ker \phi_T=V$.
Then it is not hard to show that the smallest parallelotope of the adeles $\A_K$ containing the $T$-torsion submodule $V$ has volume $q^{h((\lambda_1:\cdots:\lambda_r))}$, where $h(P)$ denotes the logarithmic Weil height of a point $P \in {\mathbb P}^{r-1}(K)$, and this height can be arbitrarily large.

We cannot apply the method of Section~\ref{proofs3} to obtain uniform bounds for torsion of Drinfeld modules of rank $r \ge 2$ until we have more information about the Galois representations arising from such Drinfeld modules.
It is perhaps reasonable to expect that the index of the image of the representation is finite for rank $r$ Drinfeld modules $\phi$ with $\End \phi = A$, but to use this to obtain uniform boundedness of torsion, it would probably be necessary to bound the index uniformly over Drinfeld modules of a given rank.
Not much is known regarding uniform boundedness of the index in the classical case either, for the image of Galois for non-CM elliptic curves.

\medskip

Let us now give the proof of Theorem~\ref{maninlike}.
Manin's proof of the number field analogue required an application of the Mordell conjecture, which had not been proved at the time, but he was able to verify the conjecture in the cases he needed.
Similarly, we will need the Mordell conjecture over characteristic $p$ function fields, which was proved by Samuel~\cite{samuel}.
I thank Ernst-Ulrich Gekeler for explaining that the necessary hypotheses for the application of the Mordell conjecture (genus at least~2, and non-isotriviality) were essentially known.

\begin{proof}[Proof of Theorem~\ref{maninlike}]
By Corollary~2.20 in~\cite{gekelerpreprint}, if $n \ge 4$, the genus of the Drinfeld modular curve $X_0(\pp^n)$ exceeds~2, and then the same will hold for the covering $X_1(\pp^n)$.
In order to apply the characteristic $p$ Mordell conjecture, we also need to know that $X_1(\pp^n)$ is non-isotrivial; i.e. that it cannot be defined over a finite field, even after base change.
Let $J$ be the Jacobian of $X_1(\pp^n)$ for some $n \ge 4$.
By~\cite{gekelerreversat}, $J$ has multiplicative reduction at $\infty$.
Then the criterion of Ogg-N\'{e}ron-Shafarevich implies that $J$ cannot be defined over a finite field, and so neither can $X_1(\pp^n)$.

By the Mordell conjecture for characteristic $p$ function fields (Theorem~4 in~\cite{samuel}), $X_1(\pp^n)$ has only finitely many $L$-rational points.
In particular, there are only finitely many $j \in L$ for which there might exist a rank~2 Drinfeld module (of that $j$-invariant) with a $\pp^n$-torsion point.
Hence it will suffice to prove that for fixed $j \in L$, there is a uniform bound for the $\pp$-primary torsion of Drinfeld modules of that $j$-invariant.
All rank~2 Drinfeld modules of a given $j$-invariant are twists of one another, so their torsion is uniformly bounded by Theorem~\ref{general}.
\end{proof}

Theorem~\ref{maninlike} could probably be generalized to general $A$.
This would involve calculating the genus of the Drinfeld modular curves for these $A$, or at least verifying that the genus is eventually at least~2.
Such a calculation can be done by starting with the genus formula for $X(1)$ (VI~5.8 in~\cite{gekelerbook}), and applying the Hurwitz formula, using the computation of the contributions of the ramified points in~\cite{gekelerthesis}.

Finally, it should be mentioned that Tamagawa~\cite{tamagawa} has developed a theory for the rank~2 case over $A=\Fq[T]$ analogous to Mazur's theory of Eisenstein ideals for modular curves over $\Q$.
In particular, he has succeeded in defining the Eisenstein quotient of the Drinfeld modular curves $X_0(\nn)$ over $K$, and has proved the finiteness of the Mordell-Weil group and the prime-to-$q$ part of the Shafarevitch-Tate group for the ``minus part'' of this Eisenstein quotient.
It is not yet clear whether this will lead to information on rational points on the Drinfeld modular curves themselves, which is what is needed to verify the uniform boundedness conjecture for rank~2 Drinfeld $\Fq[T]$-modules.
Perhaps it will be necessary to construct analogues of the ``winding quotient,'' which is used in~\cite{merel} for the classical case.

\section*{Acknowledgements}

I thank Dan Abramovich for suggesting that I try to prove Theorem~\ref{potential}.
I thank J.\ Felipe Voloch for correcting my history, and for helping me find references for the characteristic $p$ Mordell conjecture.
Finally, thanks go to Ernst-Ulrich Gekeler for explaining the other results needed in the proof of Theorem~\ref{maninlike}.


\end{document}